\documentclass[11pt]{amsproc}
\usepackage{CJK}
\usepackage{amssymb}
\usepackage{amsmath}
\usepackage{amsfonts}
\usepackage{graphicx}
\usepackage{amsmath,amstext,amsbsy,amssymb}

\newtheorem{theorem}{Theorem}[section]
\newtheorem{lemma}[theorem]{Lemma}

\newtheorem{proposition}[theorem]{Proposition}

\theoremstyle{definition}
\newtheorem{definition}[theorem]{Definition}

\theoremstyle{remark}

\numberwithin{equation}{section}
\begin{document}

\title[Generalized McKay quivers of rank three]
{Generalized McKay quivers of rank three}
\author{Xiaoli Hu}
\address{Hu: School of Sciences,
South China University of Technology, Guangzhou 510640, China}
\email{xiaolihumath@163.com}
\author{Naihuan Jing$^*$}
\address{Jing: Department of Mathematics,
   North Carolina State University,
   Raleigh, NC 27695, USA and School of Sciences, South China University of Technology, Guangzhou 510640, China}
\email{jing@math.ncsu.edu}
\author{Wuxing Cai}
\address{Cai: School of Sciences,
South China University of Technology, Guangzhou 510640, China}
\email{caiwx@scut.edu.cn}
\thanks{*Corresponding author: jing@math.ncsu.edu. Jing gratefully acknowledges the support from
NSA grant MDA 904-97-1-0062 and NSFC grant 10728102.}

\begin{abstract}
For each finite subgroup $G$ of ${\textstyle SL}_n(\mathbb C)$, we
introduce the \emph{generalized Cartan matrix} $C_{G}$ in
view of McKay correspondence from the fusion rule of
its natural representation. Using group theory, we show that the {\em generalized
Cartan matrices} have similar favorable properties such as positive
semi-definiteness as in the classical case of affine Cartan matrices (the case of
$SL_2(\mathbb C)$.
The complete McKay quivers for $SL_3(\mathbb C)$ are explicitly described and
classified based on representation theory.

\end{abstract}
\maketitle

\section{Introduction}
McKay correspondence reveals deep relations and connections
among subgroups $G$ of $SL_2(\mathbb C)$, classical Lie groups of
simply-laced types,
and resolution of singularities of $\mathbb C^2/G$ \cite{Slo, Ste}.
It has long been expected that such a
correspondence can be generalized to the special linear group
$SL_n(\mathbb C)$, and there are already several important works in
this direction \cite {IN, Ito;Nak, I, BG}. However the
complete picture is still not clear, and most results in this direction
have focused on some special subgroups such as abelian subgroups.
In this paper we give an algebraic description
of the complete generalized McKay correspondence of rank three.

One formulation is from representation theoretic
viewpoint. We take the original
elementary approach of McKay \cite{McKay} to study the fusion rule of subgroups
of $SL_n(\mathbb C)$, and in particular we focus on
the case of $n=3$. We find that to certain degree some of the
beautiful phenomena in two dimensional case are still visible in
$n$ dimensional case. The diagrams we constructed are in general some
digraphs, and the associated undirected graphs, when looking from a
distance, are similar to ADE cases.

Let $G$ be a finite nontrivial subgroup of $SL_{n}(\mathbb{C})$. We
denote the embedding of $G$ into $SL_{n}(\mathbb{C})$ as $\pi$, and
let $\rho_{0},\cdots,\rho_{s}$ be the complete list of inequivalent
irreducible representations of $G$ over $\mathbb{C}$. We then
construct a $\pi_{G}$-digraph with the irreducible representations
of $G$ as nodes and $m_{i,j}$ directed edges from $\rho_{i}$ to
$\rho_{j}$, where $\pi\otimes \rho_{i}=\bigoplus_{j}m_{ij}\rho_{j}$.
We make a convention that an undirected edge between $\rho_{i}$ and
$\rho_{j}$ represents the pair of arrows between $\rho_{i}$ and
$\rho_{j}$. This digraph is called the {\it generalized McKay
quiver} associated with $\pi$. We prove in general that all the
corresponding matrices are positive semi-definite and also a
non-trivial eigenvector is the dimension vector, which defines the
specialty of affine Cartan matrices in rank two case.

Our method can produce the generalized McKay correspondence in any dimension.
In the second part of the paper we determine the complete list of generalized McKay quivers
in rank three, including all the exceptional graphs in this case. Of course this list also includes the
 usual extended affine Dynkin diagrams
as special cases, and if one blurs some of the nodes and edges in
the quivers they also look much like the rank two cases.

We remark that a different approach to McKay correspondence in rank 3 was recently
given in \cite{But;Per} where a partial list of McKay graphs is produced. Our approach is
completely different and we have obtained the complete list of McKay quivers in rank three,
and we also indicate how to produce similar McKay graphs in higher dimensions.

One profound problem for McKay correspondence in dimension three is
to bridge the geometric approach \cite{BKR, BG} and the algebraic
method in this paper. We only discuss the algebraic picture of the
so-called McKay correspondence in rank 3 in this paper.
Nevertheless, there seems to be some partial progress between
algebra representation theory \cite{AR} and the generalized McKay
quivers, most results so far are about some subcases (such as
abelian subgroups of $SL_3(\mathbb C)$), it would be extremely
important to understand the intrinsic connection between our
digraphs and quiver theory of representation theory.

\section{The generalized McKay-Cartan matrices}
Let $G$ be a finite nontrivial subgroup of $SL_{n}(\mathbb{C})$, and
let $\gamma_i$ be its irreducible characters afforded by $\rho_i$ on
the space $V_i$ $(i=0, \cdots, |G^{*}|=s)$. The space $R(G)$ of
complex class functions on $G$ has $\{\gamma_i\}$ as a canonical
basis. For $f, g\in R(G)$, the inner product is defined by
\begin{equation}
\langle f,g\rangle=\frac{1}{|G|}\sum_{x\in
G}f(x)g(x^{-1})=\sum_{C\in G_{*}} \xi_C^{-1}f(C)g(C^{-1}),
\end{equation}
where as a convention we denote by $\xi_{C}=|Z_{G}(C)|$, and $Z_{G}(C)$
is the centralizer of one of the elements in $C$.

As usual we denote by $G_{*}$ the set of conjugacy classes: $\{C_0, C_1, $ $\cdots$, $C_{s}\}$,
 where $C_0=\{1\}$. The orthogonality
relation of irreducible characters now reads
\begin{equation}
\begin{split}
&\sum_{C}\xi_{C}^{-1}\gamma_{i}(C)\gamma_{j}(C^{-1})=\delta_{ij},\\
&\sum_{j}\gamma_{j}(C^{'})\gamma_{j}(C^{-1})=\delta_{CC^{'}}\xi_{C}.
\end{split}
\end{equation}
Let $\pi$ be the natural representation of $G$ on $\mathbb{C}^{n}$. The tensor product $\pi\otimes\rho_i$ decomposes itself into irreducible representations
\begin{equation}
\pi\otimes\rho_i\cong\bigoplus_jm_{ij}\rho_j.
\end{equation}
Then,
\begin{equation}
(nId-\pi)\otimes\rho_i\cong\bigoplus_jb_{ij}\rho_j,
\end{equation}
where $b_{ij}=n\delta_{ij}-m_{ij}.$ The matrix $M=(m_{ij})$ is
called the adjacency matrix, and $B=(b_{ij})$ the \emph{pre-Cartan
matrix} for the finite group $G$ associated to $\pi$.

\begin{definition}
The matrix $A=B+B^{t}$ is called the \emph{generalized Cartan
matrix} for the finite group $G$ associated with $\pi$.
\end{definition}
We remark that $A$ is twice of the extended Cartan matrix in the case of $n=2$.

Suppose $\pi$ is afforded by the representation $V$, the dual representation $\pi^{*}$ is the $G$-module afforded by $V^{*}=Hom_{\mathbb{C}}(V,\mathbb{C})$ and the action is given by
\begin{equation}
x.f(v)=f(x^{-1}v), x\in G,v\in V.
\end{equation}

Let $\chi_{\pi}$ be the character of $\pi$, it is well-known that
\begin{equation}
\chi_{\pi}(x)=\chi_{\pi^{*}}(x^{-1}).
\end{equation}
We also recall a well-known fact: if $U$, $V$, $W$ are $G$-modules,
then
\begin{equation}
Hom(U^{*}\otimes V, W)\cong Hom(V, U\otimes W)
\end{equation}
as $G$-modules. Hence $\langle U^{*}\otimes V,W\rangle =\langle
V,U\otimes W\rangle .$

\begin{lemma}\label{B}
$\pi^{*}\bigotimes\rho_i\cong\bigoplus_j m_{ji}\rho_j.$
\end{lemma}
Proof: Suppose $\pi^{*}\otimes\rho_i\cong\oplus_jm_{ij}^{'}\rho_j$, then
\begin{equation*}
m_{ij}^{'}=\langle \pi^{*}\otimes\rho_i, \rho_j\rangle =\langle \rho_i,
\pi\otimes\rho_j\rangle =m_{ji}.
\end{equation*}

\begin{proposition}
Let $\rho$ be a representation of $G$, and suppose
$\rho\otimes\rho_i\cong\bigoplus_jb_{ij}\rho_j.$ Let
$p_i=(\gamma_0(C_i), \cdots, \gamma_{s}(C_i))^t\in \mathbb{C}^{s}$,
then $p_i$ is an eigenvector of $B=(b_{ij})$ with eigenvalue
$\chi_{\rho}(C_{i})$.
\end{proposition}
Proof: By the orthogonality of $\gamma_i$ one has that

\begin{equation}
b_{ij}=\langle \rho\otimes\rho_i, \rho_j\rangle
=\sum_C\xi_C^{-1}\chi_{\rho}(C)\gamma_i(C)\gamma_j(C^{-1}).
\end{equation}

We compute that
\begin{equation}
\begin{split}
\sum_jb_{ij}\gamma_j(C_k)&=\sum_{j,C}\xi_C^{-1}\chi_{\rho}(C)\gamma_i(C)\gamma_j(C^{-1})\gamma_j(C_k)\\
&=\sum_C\chi_{\rho}(C)\gamma_i(C)\delta_{C,C_k}\\
&=\chi_{\rho}(C_k)\gamma_i(C_k),
\end{split}
\end{equation}
which means that $Bp_k=\chi_{\rho}(C_k)p_k.$

As a special case, the dimension vector $\delta=(dim\gamma_0,
\cdots, dim\gamma_{s})^{t}$ is an eigenvector with eigenvalue
$dim\rho$.

\begin{theorem}
The \emph{generalized Cartan matrix} $A=B+B^{t}$ is positive
semi-definite, and the dimension vector $\delta$ is an eigenvector
with eigenvalue zero.
\end{theorem}
Proof: Since $A$ is the matrix of $\rho\otimes\gamma_i=\sum
a_{ij}\gamma_j$, where $\rho=(nId-\pi)\oplus(nId-\pi^*)$, we claim
that $A$ is positive semi-definite if and only if $\langle
\rho\otimes\gamma, \gamma\rangle \geq 0$ for any $\gamma\in R(G)$.
In fact, let $\gamma=\sum x_i\gamma_i,$ $X^{T}=(x_0, \cdots,
x_{s}),$ then
\begin{equation}
\begin{split}
\langle \chi_{\rho}\otimes\gamma, \gamma\rangle
&=\sum_{i,j}x_i\langle \chi_{\rho}\otimes\gamma_i, \gamma_j\rangle
\overline{x_j}
=\sum_{i,j}x_ia_{ij}\overline{x_j}=X^{T}A\overline{X}.
\end{split}
\end{equation}

On the other hand, let $\chi$ be the any real character of degree $n$, then
\begin{equation}
\langle
\chi\otimes\gamma,\gamma\rangle=\frac{1}{|G|}\sum_{g\in
 G}\chi(g)\gamma(g)\gamma(g^{-1}) \leq n\langle
\gamma,\gamma\rangle .
\end{equation}
Hence $\langle \chi_{\rho}\otimes\gamma, \gamma\rangle = 2n\langle \gamma,
\gamma\rangle -\langle (\chi_{\pi}+\chi_{\pi^*})\otimes\gamma,\gamma\rangle
\geq0.$

\section{The $SL_3(\mathbb{C})$ McKay digraphs}

The classification of finite subgroups of  $SL_{3}(\mathbb{C})$ (up
to isomorphism) was given by Blichfelt \cite{Bli} in 1917 and completed
in \cite{YY} (see also \cite{KZSYY}). There are
$12$ types of nontrivial subgroups of $SL_{3}(\mathbb{C})$, and each type consists of
several non-isomorphic subgroups that bear similar properties. We will construct the
$\pi_G$-graphs as defined in Section 1 for all types. For each
group we first find its
irreducible representations and then decompose their tensor products with
the natural representation to construct the corresponding \emph{generalized
Cartan matrix} and the $\pi_G$-graph.
Except for the type of subgroups that are
essentially the same as the $SL_2(\mathbb C)$ case, we show each $\pi_G$-graph
together with the corresponding matrix. Some types may have
more than one graphs.

The following result will be useful in construction of irreducible
representations of semidirect products.
\begin{lemma} \label{HT}(Wigner-Mackey's Little Group method) \cite{Serre}
Suppose $G$ is the semi-direct product of a subgroup $T$ by a normal
abelian subgroup $H$. \vskip 0.1in Let $X=Hom(H,\mathbb{C}^{*})$.
Then $G$ acts on X by $(s\chi)(h)=\chi(s^{-1}hs)  ~\hbox{for}~  s\in
G,\chi\in X, h\in H$. Let $(\chi_{i})_{i\in X/T}$ be a system of
representatives for the orbits of $T$. For each $i\in X/T$, let
$T_{i}=\{t\mid t\chi_i=\chi_i\}$ and $G_i=HT_i$ be the corresponding
subgroups of $G$. Extend the function $\chi_i$ to $G_i$ by setting
$\chi_i(at)=\chi_i(a) ~\hbox{for}~ a\in H, t\in T_i$. Let $\rho$ be
an irreducible representation of $T_i$ and $\widetilde{\rho}$ be the
composition of $\rho$ with the canonical projection
$G_i\longrightarrow T_i$. Finally, let
$\theta_{i,\rho}=(\chi_i\otimes
\widetilde{\rho})\uparrow_{G_i}^{G}$. Then \vskip 0.1in

$(a)$ $\theta_{i,\rho}$ is irreducible;

$(b)$ If $\theta_{i,\rho}$ and $\theta_{i^{'},\rho^{'}}$ are
isomorphic, then $i=i^{'}$, and $\rho \simeq \rho^{'};$

$(c)$ Every irreducible representation of $G$ is isomorphic to one
 of the $\theta_{i,\rho}$.
\end{lemma}

We start with some notations. For a natural number $n$ we denote the
primitive $n$th root of unity by $\xi_n=e^{2\pi i/n}$. Then a generator of the center
of $SL_{3}(\mathbb{C})$ is given by
\begin{equation}\label{eq-W}
W=diag(\xi_3,\xi_3,\xi_3).
\end{equation}
 For positive integers $m$ and $n$, the group
 \begin{equation}
H_{m,n}=\{diag(\xi_m^k,\xi_{n}^l,\xi_{m}^{-k}\xi_{n}^{-l})|k,l\in\mathbb{Z}\}
\end{equation}
is an abelian group of order $mn$.
Now we study the 12 types one by one.

\textbf{(Type $1$) Abelian Group $H_{m,n}$}. It is well-known that
$H_{m,n}\simeq\mathbb{Z}_{m}\otimes\mathbb{Z}_{n} $. Set
$g_{k,l}=\hbox{diag}(\xi_{m}^k,\xi_{n}^l,\xi_{m}^{-k}\xi_{n}^{-l})$
$\in H_{m,n}$ and define
\begin{equation}\rho_{i,j}(g_{k,l})=
 \xi_{m}^{ik}\xi_{n}^{jl},~\mbox{for}~i=1,\cdots, m~\mbox{and}~j=1,\cdots, n.
\end{equation}
Then $\rho_{i,j}$ realize all $mn$ ($1$-dimensional) irreducible representations of $H_{m,n}$.

Let $\pi=\rho_{1, 0}\oplus \rho_{0, 1}\oplus \rho_{-1, -1}$ be the embedding of
$H_{m,n}$ into $SL_{3}(\mathbb{C})$. It follows that
\begin{equation}
\pi\otimes\rho_{i,j}=\rho_{i+1,j}\oplus\rho_{i,j+1}\oplus\rho_{i-1,j-1}.
\end{equation}
So we get the following $\pi_{H_{m,n}}$-graph.
\begin{center}
\scalebox{0.60}[0.50]{\includegraphics{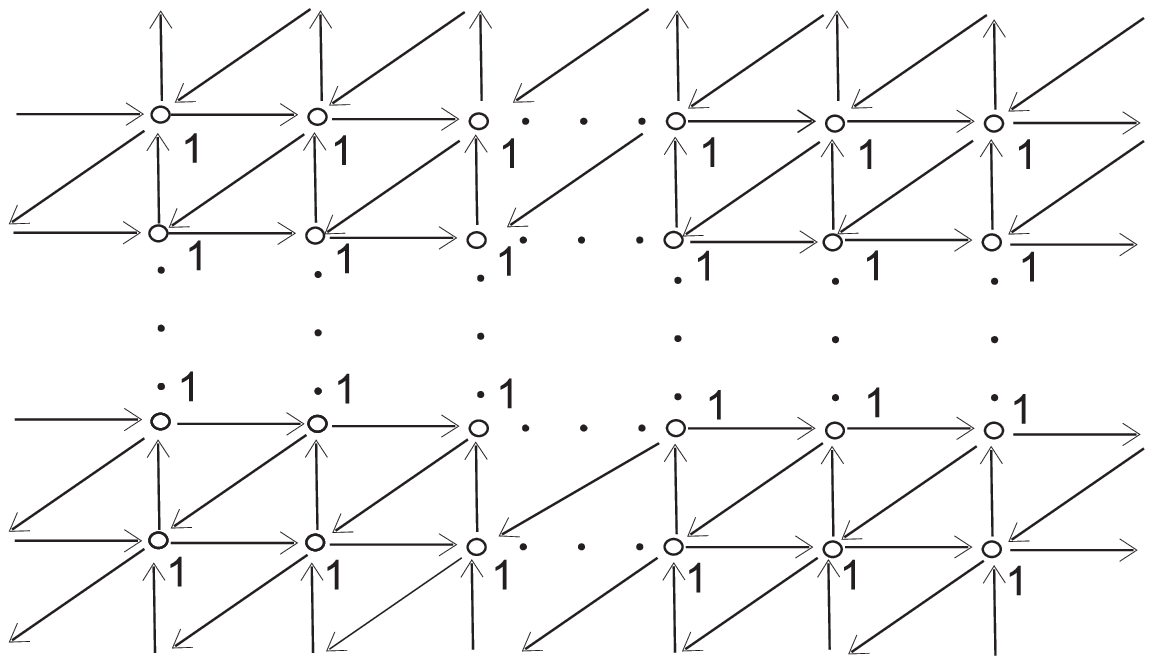}}

{\footnotesize Fig.1.}
\end{center}
Here we label each node by the dimension of the corresponding
representation. Note that the upper and lower parts of the graph are identified, similarly
the left and right parts are identified, thus we
can visualize the $\pi_{H_{m,n}}$-graph as a torus.

\begin{center}
\scalebox{0.8}[0.8]{\includegraphics{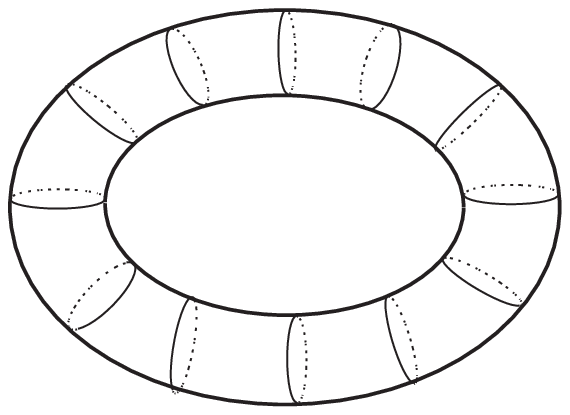}}

{\footnotesize Fig.1}
\end{center}

 We remark that in the special case of $m=1$, the graph is the
extended Dynkin diagram of $A_{n}^{(1)}$ if we erase all self
directing loops, i.e., Fig. 1-2 degenerates into a cycle.

For the \emph{generalized Cartan matrix}
$C_{H_{m,n}}=(3I-M)+(3I-M)^{T}$, where $M$ is the adjacency matrix
of $\pi_{H_{m,n}}$-graph,
 we have the following examples.
 \begin{equation}
 C_{H_{2,2}}={\left(\begin{array}{llll}
                        3  & -1 & -1& -1 \\
                        -1 & 3&-1 &-1\\
                        -1 & -1&3  &-1\\
                        -1& -1 &-1&3\\
 \end{array}\right)}.
  \end{equation}
 \begin{equation}
C_{H_{3,2}}={\left(\begin{array}{llllll}
                        6 & -2&-1 & -1 & -1& -1 \\
                        -2 & 6& -1& -1 & -1 & -1\\
                        -1 & -1& 6& -2&-1 &-1\\
                        -1 & -1&-2 & 6 &-1  &-1\\
                        -1& -1 &-1& -1& 6 & -2\\
                       -1 & -1& -1 & -1 &-2 & 6  \\
 \end{array}\right)},
 \end{equation}
 \begin{equation}
 C_{H_{3,3}}={\left(\begin{array}{lllllllll}
                        6 & -1&-1 & -1 &0& -1& 0&0&-1 \\
                        -1 & 6& -1&0& -1 & -1 & -1&-1&0\\
                        -1 & -1& 6& -1&0&-1 &0&-1&-1\\
                        -1 &0& -1& 6 &-1  &-1&-1  &-1&0\\
                        0& -1 &0& -1& 6 & -1&0&-1  &-1\\
                       -1 & -1& -1 & -1 &-1 & 6&-1&0&-1  \\
                       0&-1&0&-1&0&-1&6&-1&-1\\
                       0&-1&-1&-1&-1&0&-1&6&-1\\
                       -1&0&-1&0&-1&-1&-1&-1&6\\
 \end{array}\right)},
 \end{equation}
  \begin{equation}C_{H_{3,4}}={\left(\begin{array}{llllllllllll}
                        6 & -1& 0 & -1 & -1& -1 & 0 & 0& -1 & 0& 0& -1\\
                        -1 & 6& -1& 0 & 0 & -1& -1 & 0&-1 & -1& 0& 0\\
                        0 & -1& 6& -1& 0 & 0 & -1&-1& 0 &-1&-1& 0\\
                        -1 & 0&-1 & 6 &-1 & 0 & 0 &-1& 0 & 0&-1&-1\\
                        -1& 0 & 0 & -1& 6 & -1& 0 &-1& -1&-1& 0& 0 \\
                       -1 & -1& 0 & 0 &-1 & 6 & -1&0 & 0 &-1&-1& 0\\
                        0 &-1 &-1 & 0 & 0 &-1 & 6 &-1& 0 & 0&-1& -1\\
                        0 & 0 &-1 &-1 &-1 & 0 &-1 &6&-1 & 0& 0&-1\\
                        -1& -1& 0 & 0 &-1 & 0 & 0 &-1& 6 &-1& 0&-1\\
                         0&-1 &-1 & 0 &-1 &-1 & 0 &0 &-1 &6&-1& 0\\
                         0& 0 &-1 &-1 & 0 &-1 &-1 &0 & 0 &-1& 6&-1\\
                        -1& 0 &0  &-1 & 0 & 0 & -1&-1&-1 & 0&-1&6\\
 \end{array}\right)}.
\end{equation}\vskip 0.1in

\textbf{(Type $2$) Group ${G_{m}^{3}}$}. The group is
generated by $H_m:=H_{m,m}$ and $T$, where
\begin{equation}\label{eq-T}
T=\left(\begin{array}{lll}
                      0 &1 &0\\
                      0& 0 &1 \\
                      1 &0& 0\\
                           \end{array}\right).
\end{equation}

It is easy to see that $G_{m}^3=H_{m}\rtimes\langle T\rangle$ and
$G_{m}^3=H_{m}\uplus TH_m\uplus T^2H_m$, so we can get its
irreducible representations by Lemma \ref{HT}. Firstly, the orbits
of $X=Hom(H_{m},\mathbb{C}^{*})$ under the action of  $<T>$ are listed
as follows:
\begin{equation}
\begin{split}
\mathcal{O}(\rho_{{i,j}})=\big\{\rho_{i,j},\rho_{-j,i-j},\rho_{j-i,-i}|(i,j)\neq
(0,0),(\frac{m}{3},\frac{2m}{3}),(\frac{2m}{3},\frac{m}{3})\big\},\\
(*)\mathcal{O}(\rho_{\frac{m}{3},\frac{2m}{3}})=\{\rho_{\frac{m}{3},\frac{2m}{3}}\},
\ \
 \mathcal{O}(\rho_{\frac{2m}{3},\frac{m}{3}})=\{\rho_{\frac{2m}{3},\frac{m}{3}}\},\\
 \mathcal{O}(\rho_{0,0})=\{\rho_{0,0}\},
\end{split}
\end{equation}
where $i,j\in\{1,2\cdots,m\}$. Note: the rows labeled with $(*)$ are
excluded when m is not divisible by $3$, and the same for
(\ref{omittable4}), (\ref{omitable1}), (\ref{omitable2}) and
(\ref{omitable3}) in the following. Hence there are
$\frac{m^2-1}{3}+1$ orbits when $3\nmid m$ and $\frac{m^2-3}{3}+3$
orbits when $3|m$. Let
\begin{equation}
\begin{split}
T_1&=\left\{s\in <T>|s\cdot
\rho_{\frac{m}{3},\frac{2m}{3}}=\rho_{\frac{m}{3},\frac{2m}{3}}\right\}=<T>,\\
T_2&=\left\{s\in <T>|s\cdot
\rho_{\frac{2m}{3},\frac{m}{3}}=\rho_{\frac{2m}{3},\frac{m}{3}}\right\}=<T>,\\
T_3&=\left\{s\in <T>|s\cdot \rho_{0,0}=\rho_{0,0}\right\}=<T>,
\end{split}
\end{equation}
and let $\rho^j: T\longrightarrow \xi_3^j  (j=1,2,3)$ be the
1-dimensional irreducible representations of subgroup $<T>$.
Therefore the Little Group method gives that
\begin{align}
 \theta_{0,0,j}&=\rho_{0,0}\otimes \rho^j, \\  \label{omittable4}
 (*)\theta_{
\frac{m}{3},\frac{2m}{3},j}&=\rho_{\frac{m}{3},\frac{2m}{3}}\otimes
\rho^j, \ \
\theta_{\frac{2m}{3},\frac{m}{3},j}=\rho_{\frac{2m}{3},\frac{m}{3}}\otimes\rho^j
\end{align}
are irreducible representations of $G_m^3$ of degree 1. We can check
that

\begin{equation}
\left\{s\in<T>|s\cdot\rho_{i,j}=\rho_{i,j},(i,j)
\neq(0,0),(\frac{m}{3},\frac{2m}{3}),(\frac{2m}{3},\frac{m}{3})\right\}
=Id,
\end{equation}
then $\theta_{i,j}=(\rho_{i,j}\otimes id)\uparrow_{H_m}^{G_m^3}$ are
3-dimensional irreducible representations of $G_m^3$, where $id$ is
the trivial representation. Hence there
are nine $1$-dimensional and $\frac{m^{2}-3}{3}$ $3$-dimensional
irreducible representations in the case of $3\mid m$, and three
$1$-dimensional and $\frac{m^{2}-1}{3}$ $3$-dimensional irreducible
representations if $3\nmid m$. For $n=1, 2, 3$ and the indexes
$(i,j)\in \{1,2,\cdots,m\}\times \{1,2,\cdots,m\}\setminus
\{(0,0),(\frac{m}{3},\frac{2m}{3}),(\frac{2m}{3},\frac{m}{3})\}$,
let
\begin{equation}\label{omitable1}
\begin{split}
&\theta_{0,0,n}(T)=\xi_3^{n},\ \ \ \ \ \theta_{0,0,n}(g_{k,l})=1,\\
*~&\theta_{\frac{m}{3},\frac{2m}{3},n}(T)=\xi_3^{n},\ \ \theta_{\frac{m}{3},\frac{2m}{3},n}(g_{k,l})=\xi_3^{k+2l},\\
*~&\theta_{\frac{2m}{3},\frac{m}{3},n}(T)=\xi_3^{n},\ \ \theta_{\frac{2m}{3},\frac{m}{3},n}(g_{k,l})=\xi_3^{2k+l},\\
&\theta_{i,j}(T)=T,\qquad \
\theta_{i,j}(g_{k,l})=\hbox{diag}(\xi_m^{ik+jl},\xi_m^{(j-i)k-il},\xi_m^{-jk+(i-j)l}).
\end{split}
\end{equation}
This gives all the irreducible representations of $G_{m}^{3}$.
Recall that $\pi$, the natural representation from $G_m^3$ to
$SL_3(\mathbb{C})$, can be given in this type by
\begin{equation}\pi(g_{k,l})=\hbox{diag}(\xi_m^k,\xi_m^l,\xi_m^{-k-l}),~~~~\pi(T)=T.
\end{equation}
By computing characters, we find the following decompositions of the tensor product
of the natural representation and
the irreducible representations of $G_m^3$:
\begin{equation} \label{omitable2}
\begin{split}
&\pi\otimes\theta_{0,0,n}= \theta_{0,1,n},\ \
\pi\otimes\theta_{i,j}= \theta_{i-1,j-1}\oplus\theta_{i,j+1}
\oplus\theta_{i+1,j},\\
*~&\pi\otimes\theta_{\frac{m}{3},\frac{2m}{3},n}=
\theta_{\frac{m}{3}+1,\frac{2m}{3},n},\ \
\pi\otimes\theta_{\frac{2m}{3},\frac{m}{3},n}=
\theta_{\frac{2m}{3}+1,\frac{m}{3},n}.
\end{split}
\end{equation}

We remark that the decomposition can also be
confirmed by computing the inner product $\langle \pi\otimes \rho_i, \rho_j\rangle$.
In the following this can always be used in getting the fusion product.

The $\pi_{{G_{m}^{3}}}$-graphs are shown in Fig.2-1 for $m=6$ and
Fig.2-2 for $m=7$, one can draw the $\pi_{{G_{m}^{3}}}$-graphs for
any $m$ according to our recipe. Here we label each node
corresponding to the representation $\theta_{i, j}$ (resp. $\theta_{i, j, k}$) by its
subscript $(ij)$ (resp. $(ijk)$). (For example: $\theta_{0,0,1}$ is labeled by $(0 0 1)$.)

\begin{center}
\scalebox{0.5}[0.48]{\includegraphics{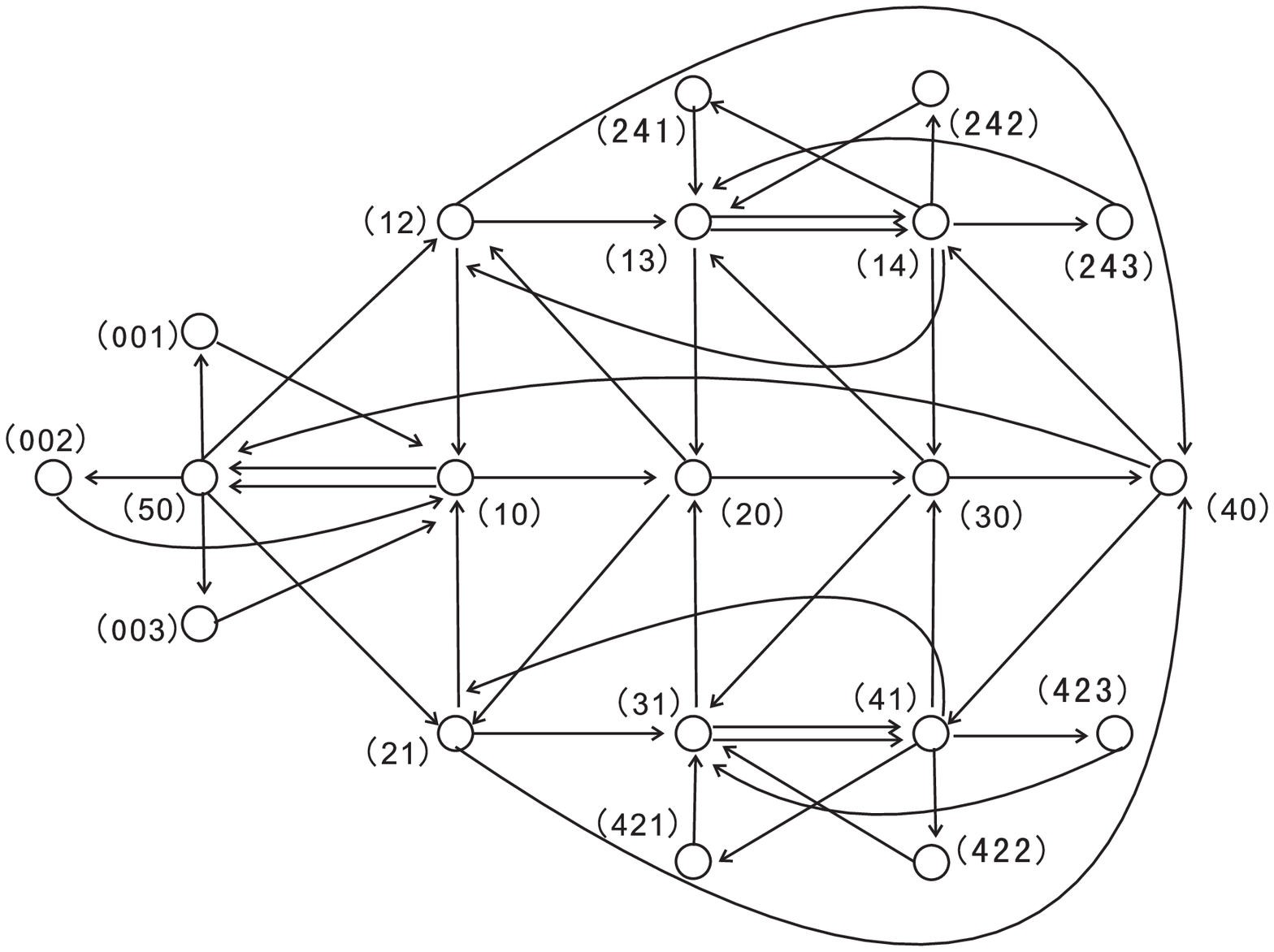}}

{\footnotesize Fig.2-1.}
\end{center}

\begin{center}
\scalebox{0.5}[0.4]{\includegraphics{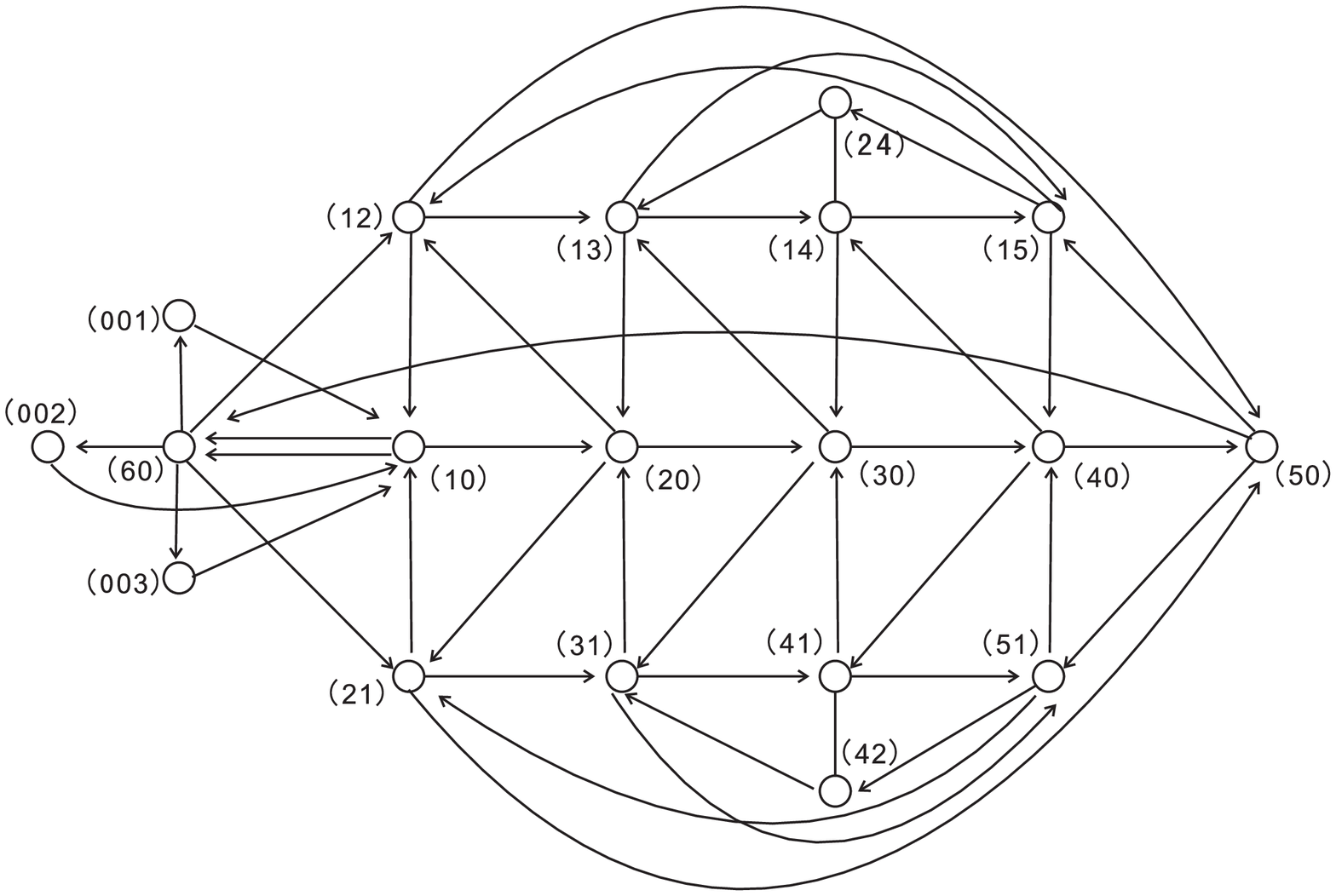}}

{\footnotesize Fig.2-2.}
\end{center}

 \textbf{(Type $3$) Group ${G_{m}^{6}}$}. The group
${G_{m}^{6}}$ is generated by $H_m$, $T$ and $R$, where $T$ was
defined in (\ref{eq-T}) and
\begin{equation}\label{eq-R}
R=\left(
\begin{array}{lll}
                      a &0 &0\\
                      0& 0 &b \\
                      0 &c& 0\\
                           \end{array}
                          \right), abc=-1.
\end{equation}

For simplicity we only consider the case of $a=b=c=-1$, which is enough to
show the main feature of this type. It is not difficult to see that
$G_{m}^6=H_{m}\rtimes S_3$, so its irreducible representations are
obtained again by using Lemma \ref{HT}. The orbits of
$X=Hom(H_m,\mathbb{C}^{*})$ under the action of $S_3$ are listed as follows:
\begin{equation}\label{omittable5}
\begin{split}
\mathcal{O}(\rho_{{i,j}})=\big\{\rho_{i,j},\rho_{j,i},
\rho_{-j,i-j},\rho_{i-j,-j},\rho_{j-i,-i},\rho_{-i,j-i}|i,j=0,\cdots,m-1\\
 \hbox{and}~(i,j)\notin \{(k,0),(0,l),(h,h),(\frac{m}{3},\frac{2m}{3}),(\frac{2m}{3},\frac{m}{3})\}\big\},\\
 \mathcal{O}(\rho_{i,0})=\{\rho_{i,0},\rho_{0,i},\rho_{-i,-i}|i\in\{1,2,\cdots,m-1\}\} ,\\
(*)~\mathcal{O}(\rho_{\frac{m}{3},\frac{2m}{3}})=\{\rho_{\frac{m}{3},\frac{2m}{3}},\rho_{\frac{2m}{3},\frac{m}{3}}\},\\
\mathcal{O}(\rho_{0,0})=\{\rho_{0,0}\},
\end{split}
\end{equation}
where $k,l,h=1,\cdots,m-1$. The irreducible
representations of $G_m^6$ can be similarly constructed as in the case of $G_m^3$, and we obtain that
 (a) If
$3|m$, there are two $1$-dimensional, four 2-dimensional, $2(m-1)$
$3$-dimensional and $\frac{m^{2}-3m}{6}$ $6$-dimensional irreducible
representations. (b) If $3\nmid m$ , there are two $1$-dimensional,
one $2$-dimensional, $2(m-1)$ $3$-dimensional and
$\frac{m^{2}-3m+2}{6}$ $6$-dimensional irreducible representations.

The irreducible representations corresponding to the orbit
$\mathcal{O}(\rho_{0,0})$:
\begin{equation}\label{omittable6}
\begin{split}
&\theta_{0,0,1}(g_{k,l})=1, \ \ \theta_{0,0, 1
}((12))=\theta_{0,0,1 }((23))=-1.\\
&\theta_{0,0,2}(g_{k,l})=1,\ \
\theta_{0,0,2}((12))=\theta_{0,0,2}((23))=1.\\
&\theta_{0,0,3}({g_{k,l}})=E_{2\times2},
\theta_{0,0,3}{(12)}=\left(\begin{matrix}
                        -1&-1\\
                          0&1\\
                           \end{matrix}
                          \right),
 \theta_{0,0,3}{(23)}=\left(\begin{array}{ll}
                         0&1\\
                          1&0\\
                           \end{array}
                          \right).\\
 \end{split}
 \end{equation}
The irreducible representations corresponding to the orbit
$\mathcal{O}(\rho_{\frac{m}{3},\frac{2m}{3}})$:
\begin{equation}
\begin{split}
&*\theta_{\frac{m}{3},\frac{2m}{3},n_1}(g_{k,l})=diag(\xi_3^{(k+2l)},\xi_3^{(2k+l)})\ \ \ (n_1=1,2,3),\\
&
*\theta_{\frac{m}{3},\frac{2m}{3},n_1}{(12)}=\left(\begin{array}{ll}
                         0&1\\
                          1&0\\
                           \end{array}
                          \right),\ \
\theta_{\frac{m}{3},\frac{2m}{3},n_1}{(23)}=\left(\begin{array}{ll}
                         0&\xi_3^{n_1}\\
                          \xi_3^{2n_1}&0\\
                           \end{array}
                          \right),
\end{split}
 \end{equation}
 The irreducible representations corresponding to the
orbit $ \mathcal{O}(\rho_{i,0})$:
\begin{equation}
\begin{split}
&\theta_{i,0,n_2}(g_{k,l})=diag(\xi_m^{ik},\xi_m^{il},\xi_m^{i(-k-l)})
\ \ \ (n_2=1,2),\\
 &\theta_{i,0,n_2}((12))=(-1)^{n_2}\left(\begin{array}{lll}
                         0&1&0\\
                          1&0&0\\
                          0&0&1\\
                           \end{array}
                          \right),
\theta_{i,0,n_2}((23))=(-1)^{n_2}\left(\begin{array}{lll}
                         1&0&0\\
                          0&0&1\\
                          0&1&0\\
                           \end{array}
                          \right).
\end{split}
 \end{equation}
The irreducible representations corresponding to the orbit $
\mathcal{O}(\rho_{i,j})$:
 \begin{equation}
\begin{split}
&\theta_{i,j}(g_{k,l})
=diag(\xi_m^{ik+jl},\xi_m^{jk+il},\xi_m^{-ik+(j-i)l)},\xi_m^{(i-j)k-jl},\xi_m^{(j-i)k-il},\xi_m^{-jk+(i-j)l}),\\
&\theta_{i,j}((12))=\left(\begin{array}{llllll}
                         0&1&0&0&0&0\\
                          1&0&0&0&0&0\\
                          0&0&0&0&0&1\\
                          0&0&0&0&1&0\\
                          0&0&0&1&0&0\\
                          0&0&1&0&0&0
                           \end{array}
                          \right),
\theta_{i,j}((23))=\left(\begin{array}{llllll}
                         0&0&0&1&0&0\\
                         0&0&0&0&0&1\\
                         0&0&0&0&1&0\\
                          1&0&0&0&0&0\\
                          0&0&1&0&0&0\\
                          0&1&0&0&0&0
                           \end{array}
                          \right).
\end{split}
 \end{equation}
Computing the characters of the tensor product of
 each irreducible representation by the natural representation $\pi=\theta_{1,0,2}$, we find that
\begin{equation}\label{omitable3}
\begin{split}
&\pi\otimes\theta_{0,0,1}=\theta_{1,0,1}, \ \
\pi\otimes\theta_{0,0,2}=\theta_{1,0,2}, \ \ \pi\otimes\theta_{0,0,3}=\theta_{1,0,1}+ \theta_{1,0,2},\\
& \pi\otimes\theta_{i,0,n_2}=\theta_{i+1,0,n_2}\oplus
\theta_{-i,1-i}\hbox{~for~} i\neq0, \ \ n_2=1,2, \\
*~&\pi\otimes\theta_{\frac{m}{3},\frac{2m}{3},1}=\theta_{\frac{m}{3}+1,\frac{2m}{3}},\
\ \ \
\pi\otimes\theta_{\frac{m}{3},\frac{2m}{3},2}=\theta_{\frac{m}{3},\frac{2m}{3}+1},\\
*~&\pi\otimes\theta_{\frac{m}{3},\frac{2m}{3},3}=\theta_{\frac{m}{3}-1,\frac{2m}{3}-1},\\
&\pi\otimes\theta_{i,j}=\theta_{i-1,j-1}\oplus\theta_{i+1,j}\oplus\theta_{i,j+1},
\hbox{~for~} {i,j}\notin\{0,\frac{m}{3},\frac{2m}{3}\}.
\end{split}
\end{equation}
These fusion rules completely determine the $\pi_{G_{m}^{6}}$-graph.
We draw the $\pi_{G_{6}^{6}}$-graph (Fig.3-1) and
$\pi_{G_{5}^{6}}$-graph (Fig.3-2) as follows.

\scalebox{0.9}[0.8]{\includegraphics{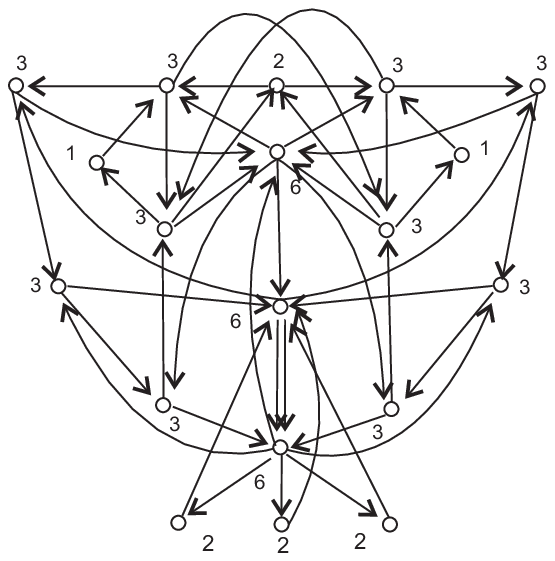}}\qquad
 \scalebox{1.0}[0.8]{\includegraphics{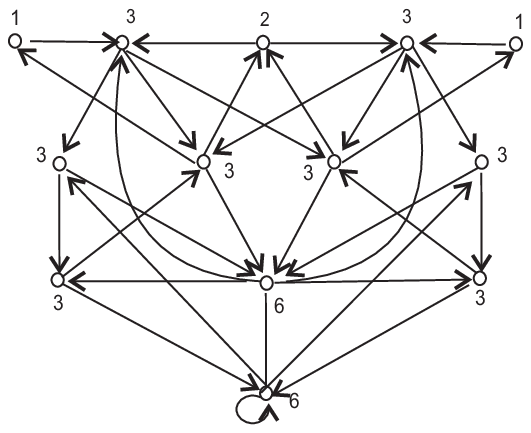}}

 \quad\qquad\qquad{\footnotesize Fig.3-1.}
\qquad\qquad\qquad\qquad \qquad\qquad{\footnotesize Fig.3-2.}

From the graph we get the adjacency matrix $M$ associated to
$\pi_{G_{m}^{6}}$-graph. Then $C(G_{m}^{6})=(3I-M)+(3I-M)^{T}$ is
the corresponding \emph{generalized Cartan matrix}. \vskip 0.1in

\textbf{(Type $4$) Group $\overline{G}_\alpha$}. For any fixed non-zero
$\alpha$ and a finite subgroup $G$ of $SL_2(\mathbb{C})$ let
$\overline{G}_\alpha=\{\hbox{diag}(\alpha^{-2},\alpha A)|A\in G\}\leq SL_3(\mathbb C)$.
 It is obvious that $\overline{G}_\alpha$ and $G$
are isomorphic. So we just need to study the finite subgroups of
$SL_{2}(\mathbb{C})$. It is well known that in the case of
$SL_{2}(\mathbb{C})$ there are five types of finite subgroups  and
the $\pi_{G}$-graphs are exactly the Coxeter graphs of
affine (simply laced) types \cite{McKay}. Correspondingly, we have five subtypes
here and the five $\overline{G}_\alpha$-graphs are also the Coxeter
diagram of affine types but with an extra circle for every
node.\vskip 0.1in

\textbf{(Type $5$) Group $G_5$}. The group $G_5$ is of order $108$
and generated by $T, S, V$, where $T$ is defined in type 2 (cf.
(\ref{eq-T})),
\begin{equation}\label{eq-S}
S=\hbox{diag}(1,\xi_{3},\xi_{3}^{2}),
\end{equation}
and
\begin{equation}\label{eq-V}
V=\frac{1}{\sqrt{-3}}\left(\begin{array}{lll}
                      1 &1 &1\\
                      1& \xi_{3}&\xi_{3}^{2} \\
                     1 &\xi_{3}^{2}&\xi_{3}\\
                           \end{array}
                          \right).
\end{equation}

    The group $G_5$ has $14$ conjugate classes. A set of conjugacy representatives is
    given by
$\{I,W,W^{2},S,ST\}\cup\{V^{i}W^{k}\mid
 1\leq i, j\leq3, 0\leq k\leq2\}, $ where
 $W$ is the generator of $Z(SL_3(\mathbb C))$ (see \ref{eq-W}).
The $14$ inequivalent irreducible representations can be described
as follows. Set $H_3=\langle T, S, W\rangle \lhd G_5$, then
$G_5=H_3\uplus VH_3 \uplus V^2H_3\uplus V^3H_3$. Let $\rho_j (j=1,
2, 3, 4)$ be the 1-dimensional irreducible representations of $G_5$:
\begin{equation}
\rho_{j}: S\longrightarrow1,\ \ T\longrightarrow1,\ \ V\longrightarrow\xi_{4}^j,\\
\end{equation}
We also need to use the canonical involution
$\sigma$ of $SL_3(\mathbb{C})$ defined by:
\begin{equation} \label{eq-sigma}
\sigma: A
\longrightarrow(A^{-1})^{T}.
\end{equation}
In terms of the natural representation $\pi$, we can construct the irreducible
representations of $G_5$ as follows. For $j = 1, 2, 3,
4$ we define
\begin{equation}
\begin{split}
\rho_{4+j}=\pi\otimes\rho_{j},\ \ \
\rho_{8+j}=\sigma\circ\rho_{4+j},\ \
\rho_{12+j}=\psi_j\uparrow_{H_3}^{G_5}.
\end{split}
\end{equation}

Again by computing their characters we obtain the following fusion rules, and the
$\pi_{G_{5}}$-graph is drawn in Fig. 5.
\begin{equation}
\begin{split}
&\pi\otimes\rho_1=\rho_5, \ \  \pi\otimes\rho_2=\rho_7,\ \
\pi\otimes\rho_3=\rho_9, \ \  \pi\otimes\rho_4=\rho_{11},\\
&\pi\otimes\rho_5=\rho_6\oplus\rho_8\oplus\rho_{10},\ \ \pi\otimes\rho_6=\rho_1\oplus\rho_{13}\oplus\rho_{14},\\
&\pi\otimes\rho_7=\rho_6\oplus\rho_8\oplus\rho_{12},\ \ \pi\otimes\rho_8=\rho_2\oplus\rho_{13}\oplus\rho_{14},\\
&\pi\otimes\rho_9=\rho_6\oplus\rho_{10}\oplus\rho_{12},\ \
\pi\otimes\rho_{10}=\rho_3\oplus\rho_{13}\oplus\rho_{14},\\
&\pi\otimes\rho_{11}=\rho_8\oplus\rho_{10}\oplus\rho_{12},\ \
\pi\otimes\rho_{12}=\rho_4\oplus\rho_{13}\oplus\rho_{14},\\
&\pi\otimes\rho_{13}=\rho_5\oplus\rho_7\oplus\rho_9\oplus\rho_{11},\\
&\pi\otimes\rho_{14}=\rho_5\oplus\rho_7\oplus\rho_9\oplus\rho_{11}.
\end{split}
\end{equation}

Similarly, we let $M$ be the adjacency matrix associated with
$\pi_{G_{5}}$-graph, and the \emph{generalized Cartan matrix}
$C(G)=(3I-M)+(3I-M)^{T}$.
\begin{center}
\scalebox{0.8}[0.8] {\includegraphics{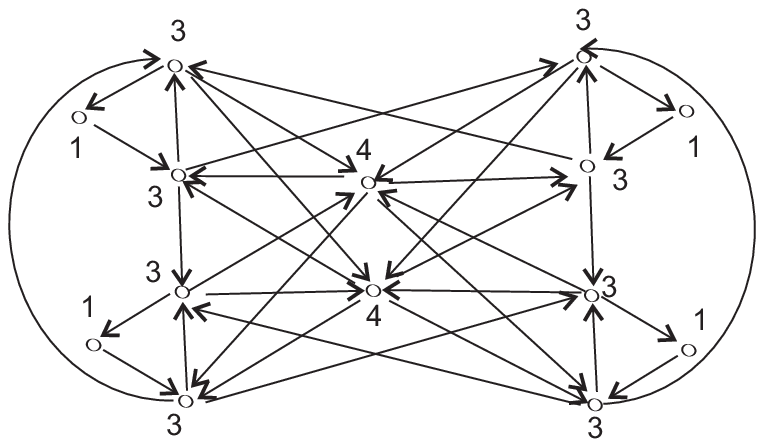}}

{\footnotesize Fig. 5}
\end{center}

\textbf{(Type $6$) Group $G_6$}. The group of order 216 is generated
by $G_5$ and $K$, where
\begin{equation}\label{eq-K}
K=\frac{1}{\sqrt{-3}}\left(\begin{array}{lll}
                      1 &1 &\xi_{3}^{2}\\
                      1& \xi_{3}&\xi_{3} \\
                    \xi_{3}&1&\xi_{3}\\
                           \end{array}
                          \right).
\end{equation}

It is known that there are $16$ conjugacy classes in $G_6$, hence
$G_6$ has $16$ irreducible representations. Since $G_5$ is a normal
subgroup of $G_6$, we can find some irreducible representations of
$G_6$ by induction. In fact we find two $6$-dimensional irreducible
representations denoted as $\rho^{14}$ and $\rho^{15}$, one
$2$-dimensional representation denoted as $\rho^{13}$ and one
$8$-dimensional irreducible representations denoted as $\rho^{16}$.
The other four $1$-dimensional and eight $3$-dimensional
irreducible representations are given as follows: for $i,j=0,1$,
\begin{equation}
\begin{split}
&\rho^{(i,j)}:T\longrightarrow 1,\quad S\longrightarrow1,\quad  V\longrightarrow \xi_2^{i},\quad  K\longrightarrow\xi_2^{j},\\
&\rho^{(2+i,2+j)}=\pi\otimes \rho^{(i,j)},\ \ \ \ \rho^{(4+i,4+j)}=\sigma\circ\rho^{(2+i,2+j)}.\\
\end{split}
\end{equation}

The natural representation is seen to be $\pi=\rho^{(2,2)}$, and we have the
following decompositions by comparing characters.
The $\pi_{G_6}-$graph is drawn in Fig.6., and the corresponding \emph{generalized Cartan matrix} is then determined accordingly.

\begin{equation}
\begin{split}
&\pi\otimes\rho^{(i,j)}=\rho^{(2+i,2+j)}, \\
&\pi\otimes\rho^{(2+i,2+j)}=\rho^{(5-i,5-j)}\oplus\rho^{14},\\
&\pi\otimes\rho^{(4+i,4+j)}=\rho^{(1-i,1-j)}\oplus\rho^{16},\\
&\pi\otimes\rho^{13}=\rho^{15}, \ \ \pi\otimes\rho^{14}=2\rho^{16}\oplus\rho^{13},\\
&\pi\otimes\rho^{15}=\rho^{9}\oplus\rho^{10}\oplus\rho^{11}\oplus\rho^{12}\oplus\rho^{14},\\
&\pi\otimes\rho^{16}=\rho^{5}\oplus\rho^{6}\oplus\rho^{7}\oplus\rho^{8}\oplus2\rho^{15}.
\end{split}
\end{equation}
\begin{center}
\scalebox{0.75}[0.75] {\includegraphics{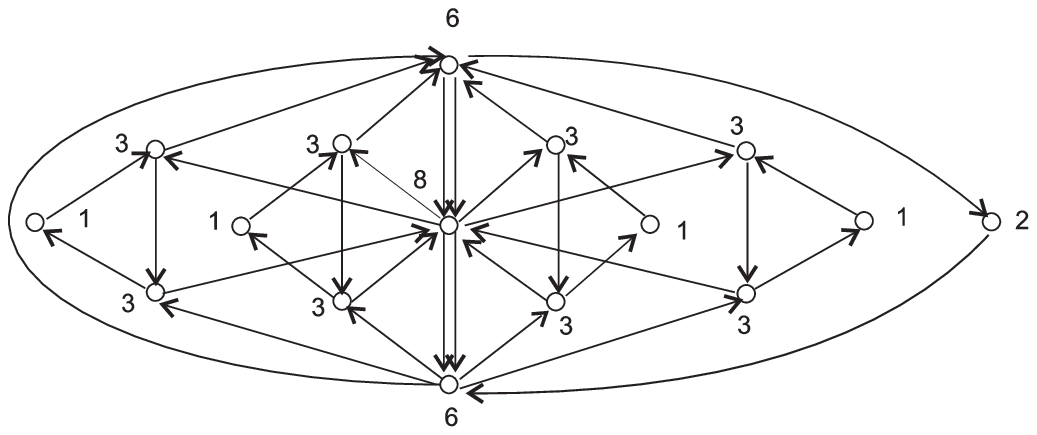}}

{\footnotesize Fig.6.}
\end{center}\vskip 0.1in

\textbf{(Type $7$) Group $G_7$}. The group $G_7$ of order $60$ is
generated by $T$, $E_{2}$, and $E_{3}$, where $T$ was defined as above
(\ref{eq-T}) and
\begin{equation}\label{eq-E}
E_{2}=\hbox{diag}(1,-1,-1), \ \ E_{3}=\frac{1}{2}\left(\begin{array}{lll}
                      -1 &\mu_{-} &\mu_{+}\\
                      \mu_{-}& \mu_{+} &-1 \\
                     \mu_{+}&-1&\mu_{-}\\
                           \end{array}
                          \right),\mu_{\pm}=\frac{1}{2}(-1\pm\sqrt{5}).
\end{equation}

The group $G_7$ is isomorphic to the simple group $A_{5}$, whose
irreducible characters are known. Table 1 lists its character table,
where $\nu_{\pm}=\frac{-1\pm\sqrt{5}}2$. For the natural character
$\chi_{\pi}=\chi_{3}$ of the imbedding, we have the following
decompositions, and the $\pi_{G_{7}}$-graph is shown in Fig.7.

\begin{table}
\begin{tabular}{|c|c|c|c|l|c|}
\hline  Class & $1$& $2$ & $3$ & $5$ &$5$\\
\hline $\chi_1$ &1  &  1  & 1 & 1 & 1  \\
\hline $\chi_2$ & 3 &   -1 & 0 & $\nu_{+}$ & $\nu_{-}$ \\
\hline $\chi_3$ &3 &-1 & 0 & $\nu_{-}$ & $\nu_{+}$\\
\hline $\chi_4$ &4 & 0 &1 & -1 &-1 \\
\hline $\chi_5$ & 5& 1 &-1 &0 &0\\
\hline
\end{tabular}
\caption{The characters of $G_7$\label{tab}}
\end{table}

\begin{table}
\begin{tabular}{|c|c|c|c|c|l|c|}
\hline  Class & $1$& $2$ & $4$ & $3$ &$7_1$&$7_2$\\
\hline $\chi_1$ &1  &  1  & 1 & 1 & 1 & 1 \\
\hline $\chi_2$ & 6 & 2 & 0 & $0$ & $-1$& $-1$ \\
\hline $\chi_3$ &7 &-1 & -1 & $1$ & $0$& $0$\\
\hline $\chi_4$ &8 & 0 &0 & -1 &1&1  \\
\hline $\chi_5$ & 3& -1 &1 &0 &$\alpha_{+}$&$\alpha_{-}$\\
\hline $\chi_6$ & 3& -1 &1 &0 &$\alpha_{-}$&$\alpha_{+}$\\
\hline
\end{tabular}
\caption{The characters of $G_8$\label{tab}}
\end{table}

\begin{equation}
\begin{split}
 \chi_{\pi}\otimes\chi_{1}&=\chi_{3},\\
 \chi_{\pi}\otimes\chi_{2}&=\chi_{2}\oplus\chi_{4}\oplus\chi_{5},\\
 \chi_{\pi}\otimes\chi_{3}&=\chi_{1}\oplus\chi_{3}\oplus\chi_{4},\\
 \chi_{\pi}\otimes\chi_{4}&=\chi_{2}\oplus\chi_{3}\oplus\chi_{4}\oplus\chi_{5},\\
 \chi_{\pi}\otimes\chi_{5}&=\chi_{2}\oplus\chi_{4},
\end{split}
C_{G_{7}}={\left(\begin{array}{lllll}
                       3&0&-1&0&0\\
                       0&2&0&-1&-1\\
                        -1&0&2&-1&0\\
                         0&-1&-1&2&-1\\
                         0&-1&0&-1&3\\
                          \end{array}\right)}.
\end{equation}
\begin{center}\scalebox{0.8}[0.75] {\includegraphics{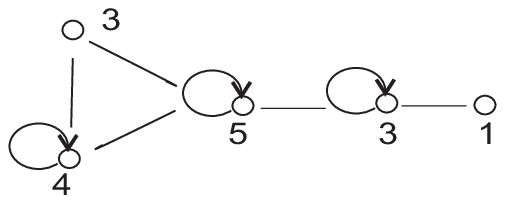}}

{\footnotesize Fig.7.}
\end{center}

\textbf{(Type $8$) Group $G_{8}$}. The group $G_8$ is the Klein simple group of order
$168$ generated by $T$ (defined in \ref{eq-T}), $X_7$ and $U$ defined below:
\begin{equation}\label{eq-XU}
X_7=\hbox{diag}(\beta,\beta^2,\beta^4),\ \
U=\frac{1}{\sqrt{-7}}\left(\begin{array}{lll}
                      a&b&c\\
                       b&c&a\\
                        c&a&b\\
                          \end{array}\right),
\end{equation}
where for $\beta=\xi_7,a=\beta^4-\beta^3,b=\beta^2-\beta^5,$ and
$c=\beta-\beta^6$.
It is known that $G_8$ has 6 irreducible
characters, whose character table is given in Table 2, where
$\alpha_{\pm}=\frac{-1\pm\sqrt{-7}}2$. The natural character
$\chi_{\pi}=\chi_5$, then we have the following decompositions of
$\chi_{\pi}\otimes \chi_{i}$$(i=1, 2,\cdots, 6)$. The
$\pi_{G_{8}}$-graph is drawn in Fig.8 and the \emph{generalized
Cartan matrix} is $C_{G_{8}}.$
\begin{equation}
\begin{array}{ll} \chi_{\pi}\otimes\chi_1=\chi_5,\\
 \chi_{\pi}\otimes\chi_2=\chi_3\oplus\chi_4\oplus\chi_6,\\
 \chi_{\pi}\otimes\chi_3=\chi_2\oplus\chi_3\oplus\chi_4,\\
 \chi_{\pi}\otimes\chi_4=\chi_2\oplus\chi_3\oplus\chi_4\oplus\chi_5,\\
 \chi_{\pi}\otimes\chi_5=\chi_2\oplus\chi_6,\\
 \chi_{\pi}\otimes\chi_6=\chi_1\oplus\chi_4,
\end{array}
C_{G_{8}}=\left(\begin{array}{llllll}
                       6&0&0&0&-1&-1\\
                       0&6&-2&-2&-1&-1\\
                        0&-2&4&-2&0&0\\
                         0&-2&-2&4&-1&-1\\
                          -1&-1&0&-1&6&-1\\
                          -1&-1&0&-1&-1&6\\
                  \end{array}\right).
\end{equation}
\begin{center}
 \scalebox{0.75}[0.75] {\includegraphics{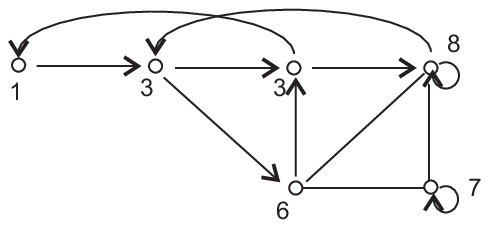}}
{\footnotesize Fig.8.}
\end{center}
In this case we can obtain the fusion rule easily by observing the character table.
For example, the product of the characters $\chi_5=\chi_{\pi}$ and $\chi_2$ is the sum
of the character (row) vectors $\chi_3$, $\chi_4$ and $\chi_6$ in Table 2.

\textbf{(Type $9$) Group $G_{9}$}. The group $G_{9}$ of order $180$
is generated by group $G_7$ and matrix $W$ (cf \ref{eq-W}). We have $G_{9}\simeq
A_5\times\langle W\rangle $ , so it is easy to get all the
irreducible representations of $G_{9}$ from those of $A_5$ and group
$\langle W\rangle$. Then the $\pi_{G_{9}}$-graph is Fig.9. and the
\emph{generalized Cartan matrix} is $C_{G_{{9}}}$.
\begin{equation}
C_{G_{9}}= \left(\begin{array}{lll}
                       6E&-B&-B\\
                       -B^{t}&6E&-B\\
                       -B&-B^{t}&6E\\
                  \end{array}\right),~
B=\left(\begin{array}{llllll}
                       0&0&1&0&0 \\
                        0&0&0&1&1 \\
                         1&0&1&0&1\\
                        0&1&0&1&1\\
                         0&1&1&1&1\\
\end{array}\right).
\end{equation}
\begin{center}
\scalebox{0.75}[0.7] {\includegraphics{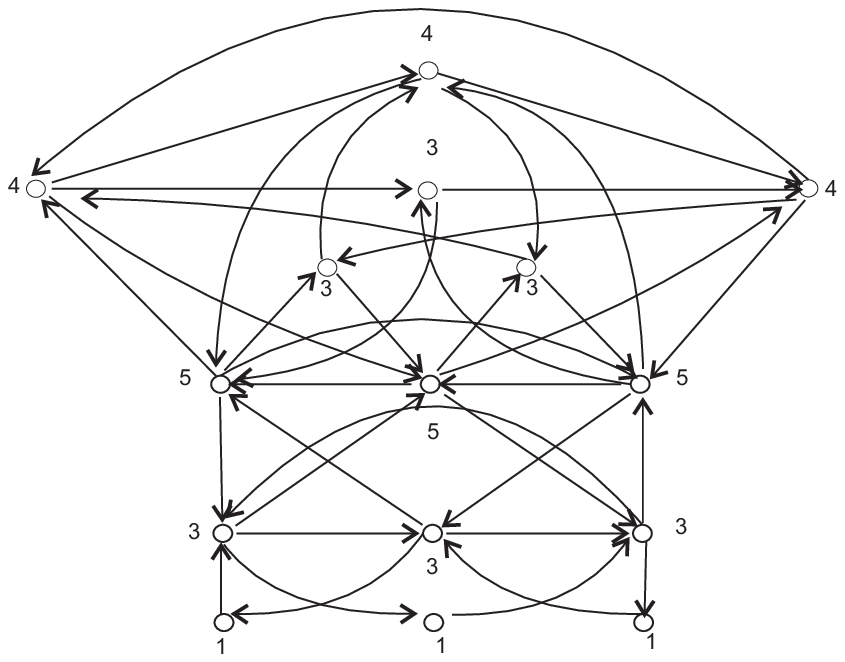}}

{\footnotesize Fig.9.}
\end{center}
\vskip 0.1in

\textbf{(Type $10$) Group $G_{10}$}. The group $G_{10}$ of order
$504$ is generated by $G_{8}$ and matrix $W$. It is easy to see that $G_{10}\simeq
G_{8}\times\langle W\rangle $. The $\pi_{G_{10}}$-graph is drawn in Fig.10
and the \emph{generalized Cartan matrix} is $C_{G_{10}}.$
\begin{equation}C_{G_{10}}= \left(\begin{array}{lll}
                       6E&-B&-B^{t}\\
                       -B^{t}&6E&-B\\
                       -B&-B^{t}&6E\\
                  \end{array}\right),~B=\left(\begin{array}{llllll}
                       0&0&0&0&0&1\\
                       0&0&1&1&1&0\\
                        0&1&1&1&0&0\\
                         0&1&1&1&0&1\\
                          1&0&0&1&1&0\\
                          0&1&0&0&1&0\\
                  \end{array}\right).
\end{equation}
\begin{center}
\scalebox{0.75}[0.7] {\includegraphics{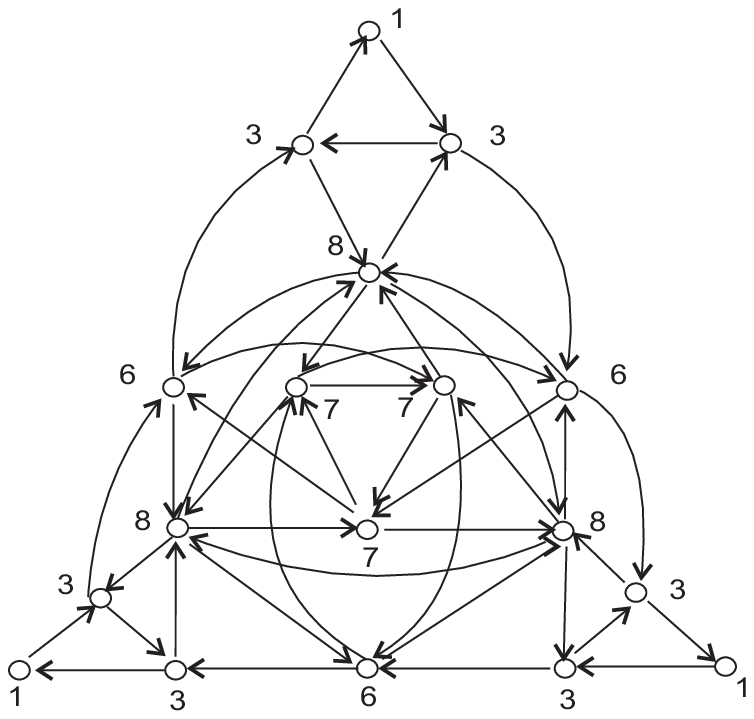}}

{\footnotesize Fig.10.}
\end{center}

\textbf{(Type $11$) Group $G_{11}$}. The $G_{11}$ of order $648$ is
generated by $G_5$ and the matrix $M$, where
\begin{equation}\label{eq-M}
M=\hbox{diag}(\varepsilon,\varepsilon,\varepsilon\xi_3),
~\varepsilon^3=\xi_3^2.
\end{equation}

 We first obtain three $1$-dimensional representations
$\rho_{i}$ by the defining relations of $G_{11}$($i=0,1,2$). We
then have three $3$-dimensional irreducible representations
$\rho_{2+i}=\pi\otimes \rho_{i}$, where $\pi$ is the natural representation.
Furthermore, we get three more
non-isomorphic $3$-dimensional irreducible representations
$\rho_{5+i}=\sigma\cdot\rho_{2+i}$, where $\sigma$ is the action of transpose inverse
(\ref{eq-sigma}). Obviously, $G_6$ is a normal subgroup of $G_{11}$,
so we can get eight $9$-dimensional representations by inducing
the $3$-dimensional irreducible representations from $G_6$ to
$G_{11}$. By Frobenius reciprocity we observe that among these induced representations, there are just
two new irreducible representations which are $9$-dimensional, and the
rest can be decomposed into direct sums of one $3$-dimensional
and one new $6$-dimensional representations. In short, there are three
$1$-dimensional, three $2$-dimensional, seven $3$-dimensional, six
$6$-dimensional, three $8$-dimensional and two $9$-dimensional
irreducible representations. Finally, we decompose the tensor
product of the natural representation $\pi$ with the irreducible
representations. The $\pi_{G_{11}}-$graph is shown in
Fig.11.
\begin{center}
\scalebox{0.8}[0.7]{\includegraphics{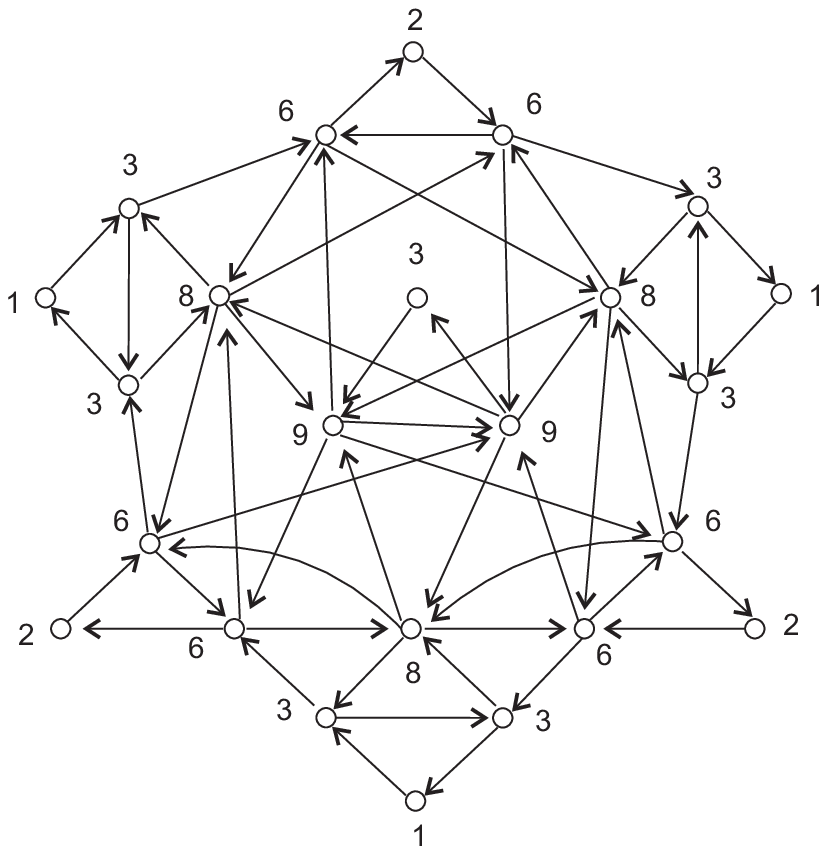}}

{\footnotesize Fig.11.}
\end{center}

\textbf{(Type $12$) Group $G_{12}$}. The group of order $1080$ is
generated by $G_7$ and matrix $E_4$, where
\begin{equation}\label{eq-E_4}
E_4=\left(\begin{array}{lll}
                       1&0&0\\
                       0&0&-w\\
                        0&-w^2&0\\
                  \end{array}\right).
\end{equation}

The center of group $G_{12}$ is $Z(G)=\langle W\rangle $ (cf. \ref{eq-W}) , and
$G_{12}/Z(G)\simeq A_6$. Thus for an irreducible representation
$\rho$ of $A_6$, we can lift it to an irreducible representation of $G_{12}$
through the canonical map. In this way we can find two
$5$-dimensional, one $9$-dimensional, one $10$-dimensional and two
$8$-dimensional irreducible representations of $G_{12}$. Moreover,
there are one $1$-dimensional and four $3$-dimensional irreducible
representations. To sum up, $G_{12}$ has one $1$-dimensional, four
$3$-dimensional, two $5$-dimensional, two $6$-dimensional, two
$8$-dimensional, three $9$-dimensional, one $10$-dimensional and two
$15$-dimensional irreducible representations. Then we decompose the
tensor product of the natural representation $\pi$ with the
irreducible representations, and the McKay quiver $\pi_{G_{12}}$ is
drawn in Fig.12.
\begin{center}
\scalebox{0.85}[0.65] {\includegraphics{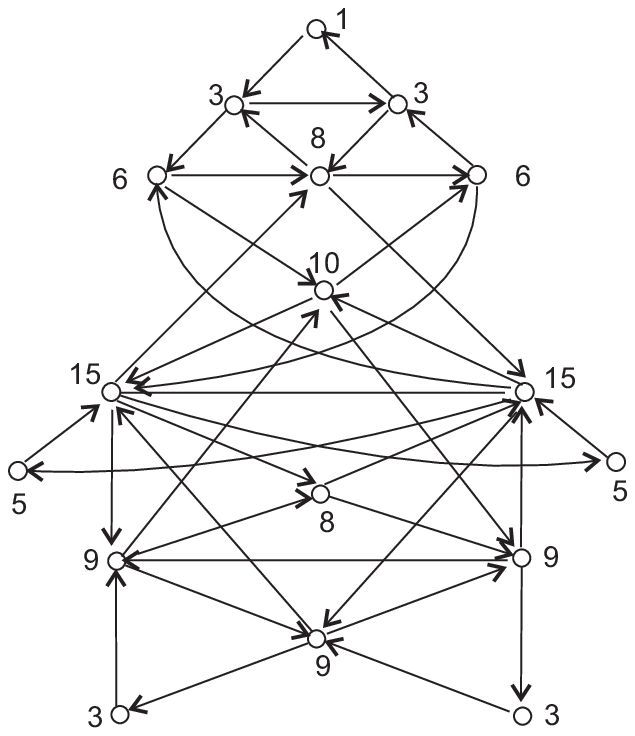}}

{\footnotesize Fig.12.}
\end{center}

\vskip 0.1in
\bibliographystyle{amsalpha}

\end{document}